\documentclass[12pt]{amsart}
\usepackage{amsthm,amsmath,amssymb,amscd,graphics,enumerate}
\usepackage{verbatim}
\usepackage[all]{xy}

{
   \newtheorem{theorem}[subsubsection]{Theorem}
   \newtheorem{proposition}[subsubsection]{Proposition}     
   \newtheorem{lemma}[subsubsection]{Lemma}

}
{\theoremstyle{definition}

   \newtheorem{definition}[subsubsection]{Definition}
   
}

\newcommand{\RR}{{\mathbb{R}}}

\newcommand{\CC}{{\mathbb{C}}}
\newcommand{\QQ}{{\mathbb{Q}}}

\newcommand{\ZZ}{{\mathbb{Z}}}

\newcommand{\GG}{{\mathbb{G}}}

\newcommand{\bR}{{\mathbf{R}}}

\newcommand{\bL}{{\mathbf{L}}}
\newcommand{\bmu}{{\boldsymbol{\mu}}}

\newcommand{\cD}{{\mathcal D}}
\newcommand{\cE}{{\mathcal E}}
\newcommand{\cF}{{\mathcal F}}
\newcommand{\cG}{{\mathcal G}}

\newcommand{\cM}{{\mathcal M}}

\newcommand{\cO}{{\mathcal O}}
\newcommand{\cP}{{\mathcal P}}

\newcommand{\cS}{{\mathcal S}}

\newcommand{\cW}{{\mathcal W}}
\newcommand{\cX}{{\mathcal X}}
\newcommand{\cY}{{\mathcal Y}}

\newcommand{\wB}{{B \otimes \omega}}
\newcommand{\wE}{{E \otimes \omega}}

\newcommand{\Spec}{\operatorname{Spec}}

\newcommand{\Per}{{\operatorname{Per}}}
\newcommand{\Hom}{{\operatorname{Hom}}}

\newcommand{\Ext}{{\operatorname{Ext}}}

\newcommand{\Tor}{{\operatorname{Tor}}}

\newcommand{\codim}{\operatorname{codim}}

\newcommand{\lrar}{\longrightarrow}

\newcommand{\Zn}{{\ZZ/n\ZZ}}

\begin{document}
\title{Flops, flips and perverse point sheaves on threefold stacks}  
\author[D. Abramovich]{Dan Abramovich}
\thanks{Research of D.A. partially supported by NSF grant DMS-0070970}  
\address{Department of Mathematics\\ Boston University\\ 111 Cummington
         \\ Boston, MA 02215\\ U.S.A.} 
\curraddr{Department of Mathematics,
Box 1917, 
Brown University, Providence, RI, 02912}
\email{abrmovic@math.bu.edu,abrmovic@math.brown.edu}
\author[J.\ C.\ Chen]{Jiun C.\ Chen}
\address{Department of Mathematics\\ Harvard  University\\ 1 Oxford
         \\ Cambridge, MA 02138\\ U.S.A.} 
\email{jcchen@math.harvard.edu}
\maketitle
\section{Introduction}

We work with varieties over $\CC$.
\subsection{Threefold flops as Moduli of perverse point sheaves}
In \cite{br}, T.\ Bridgeland considered a morphism $X\to Y$ of complex
projective varieties satisfying 
\begin{enumerate}
\item[\bf (B.1) ] $\RR f_* \cO_{X} = \cO_Y$, and
\item[\bf (B.2) ] $\dim f^{-1} \{z\} \leq 1$ for every $z\in Y$.
\end{enumerate}

Note that our notation differs from \cite{br}, in that $X$ and $Y$ are
switched. 
Bridgeland defined an abelian subcategory $Per(X/Y)\subset D^b(X)$ of the
derived category $D^b(X):= D^b_{coh}(Qcoh X)$ of bounded complexes with coherent
cohomology, and in this subcategory he identifies certain objects called {\em
perverse point sheaves}.  He proved: 
\begin{theorem}[Bridgeland, \cite{br}, Theorem 3.8] There exists a fine moduli
space $M(X/Y)$ of perverse point sheaves, which is projective over $Y$. It
contains a distinguished component $W\subset M(X/Y)$ which is birational to
$Y$. 
\end{theorem}
Bridgeland further proved the following remarkable result:
\begin{theorem}[Bridgeland, \cite{br}, Theorem 1.1]\label{Th:bridgelan-flop}
Assume that $X$ is a smooth 
threefold and $X \to Y$ is a flopping contraction. Then
\begin{enumerate}
\item $W$ is smooth,
\item The Fourier-Mukai type transformation induced by the universal perverse
point sheaf is an isomorphism, and
\item $W \simeq X^+$, the flop of $X \to Y$.
\end{enumerate}
\end{theorem}

In \cite{Chen}, Theorem \ref{Th:bridgelan-flop} is generalized to the
case where $X$ is a threefold with Gorenstein terminal singularities and $X \to
Y$ is a flopping contraction.

\subsection{Non-Gorenstein threefolds as stacks}
In this paper we are concerned with generalizing the results of  \cite{br},
\cite{Chen} to some $\QQ$-Gorenstein threefolds, using algebraic stacks, as
commented in \cite{Chen}, 
Section 1.8. It should be pointed out that Kawamata obtained some very general
results \cite{Kawamata-francia,Kawamata-D}, also using algebraic stacks,
concentrating on  equivalences 
(or embeddings) of derived categories in birational transformations. To
avoid excessive overlap with Kawamata's work, we emphasize the moduli
construction of birational transformations. 

The underlying idea is the following: in all considerations of Fourier--Mukai
transforms, smoothness is an essential assumption. Thus, if one is to prove
results analogous to \ref{Th:bridgelan-flop} for singular varieties,
some ``hidden smoothness'' would 
be desirable. In the terminal Gorensten case, a deformation $\cX \to \cY$ of $X
\to Y$ with  $\cX$ smooth is used. In the $\QQ$-Gorenstein case, the
singularities can be simplified back to the Gorenstein case, by
taking the canonical Gorenstein stack (this is also the main idea behind
Kawamata's result \cite{Kawamata-francia}). 

\subsection{Threefold terminal flops}
By way of comparison, we consider two very different cases here. The first is
that of  a threefold terminal flopping contraction $X \to Y$, with flop $X^+
\to Y$. In this case 
Kawamata  \cite{Kawamata-D}, Theorem 6.5, has proven an equivalence of
 derived categories of the canonical covering stacks $\cX\to X$
and $\cX^+ \to X^+$. Here we add a footnote to Kawamata's result: the
flop can in some sense be constructed {\em a priori} as a moduli space. 
Indeed the entire program of \cite{bkr, br,  
Chen} works:

\begin{proposition}\label{Prop:flop}
 Let $X\to Y$ be a flopping contraction of terminal threefolds,
and let $\cX \to \cY$ be the contraction of associated canonical covering
stacks. Then 
\begin{enumerate} 
\item the distinguished component $\cW(\cX/\cY)$ of the moduli stack
of perverse point sheaves has terminal Gorenstein singularities,
\item the Fourier--Mukai type transform $$D^b(\cW) \to D^b(\cX)$$ given by the
universal object is an equivalence, and
\item  $\cW(\cX/\cY)\simeq \cX^+$, the flop of $\cX \to \cY$.
\end{enumerate}
\end{proposition}
 This in particular gives a slightly different
approach to \cite{Kawamata-D}, Theorem 6.5.

While our result extends the ideas of \cite{br} to this
case, there is still an unsatisfactory point which we do not know how
to resolve: our moduli construction relies on a presentation of the
canonical covering stack; it would be desirable to have a construction
directly in terms of the stack.

\subsection{Threefold Francia flips}
The other case we consider is that of the simplest sequence of flips in
dimension 3, the so-called Francia flips, obtained as the quotient of the
standard threefold flop by a particular action of a cyclic group. In this case
something new and very different happens: consider the canonical covering stack
$\cX$ of $X$; in this case there is no canonical covering stack of $Y$ since
$Y$ is not $\QQ$-gorenstein. The usual moduli space of perverse
point sheaves $M(\cX/Y)$ is not isomorphic to $X^+$. Instead we construct a new
abelian category $Per^{(-1,0)}(\cX/Y)$ and consider the analogous notion of
perverse point sheaves. We then show:

\begin{theorem}\label{Th:francia-W} 
Let $X\to Y$ be a flipping contraction which is locally
analytically of the Francia type. Let $\cX\to X$ be the canonical covering
stack. Then there is a separated moduli space $M(\cX/Y)$ of $(-1,0)$-perverse
point sheaves on $\cX/Y$, whose distinguished component  $W(\cX/Y)$ is
projective and 
birational over $Y$.
\end{theorem}   
\begin{theorem} \label{Th:francia-FM}
\begin{enumerate} 
\item the distinguished component $W(\cX/Y)$  is smooth,
\item the Fourier--Mukai type transform $$D^b(W) \to D^b(\cX)$$ given by the
universal object is fully faithful, and
\item  $W(\cX/Y)\simeq X^+$, the flip of $\cX \to Y$.
\end{enumerate}
\end{theorem}

It should be again pointed out that the existence of a fully faithful embedding
$D^b(W) \to D^b(\cX)$ is a special case of the main theorem of
\cite{Kawamata-francia}.

Unlike the previous cases, we do not provide an {\em a priori}
construction of the flip as a moduli space: the present proof of these theorems 
uses the existence 
of $X^+$ in a very explicit way, in two points. First, our proof of projectivity of
$X/Y$ relies on the existence of a morphism $X^+ \to W$, i.e. the existence of
a family of  $(-1,0)$-perverse point sheaves for $\cX/Y$ parametrized by
$X^+$. Second, in our  proof  of the fully faithful embedding we used the same
morphism to show that $W\to Y$ does not contract a surface to a point. This is
needed when applying the method of \cite{bkr} with the intersection theorem.

\subsection{Some generalities}

We use the notation $D(X)$ exclusively for the derived category of
{\em coherent} sheaves on $X$, and $D^b(X)$ for the bounded category.

For the definition and properties of perverse sheaves, perverse ideal
and point sheaves, and their properties we rely on \cite{br} and
\cite{Chen}. We cannot improve here much on the presentation there.

For a Deligne--Mumford stack $\cX$ with coarse moduli space $\pi:\cX
\to X$ we note that $\pi_*\cO_\cX = \cO_X$, and that $\pi_*$ is exact
on quasi-coherent sheaves.
  
If $\cX = [V/G]$, with $V$ a variety and $G$ a finite group, then
$Coh(\cX) \simeq Coh^G(V)$. The moduli space is simply $X =
V/G$. Write $q: V \to X$ for the schematic quotient morphism.  If
$\cF$ is a sheaf on $\cX$ 
corresponding to a $G$-sheaf $\cG$ on $V$, then $\pi_* \cF = (q_*\cG)^G$. 

\subsection{Acknowldgements} We heartily thank Tom Bridgeland for his
help, especially with the crucial Section \ref{Sec:thankbr}.

\section{Canonical covering stacks and the moduli space for $\QQ$-Gorenstein
flops} 
\subsection{The canonical covering stack}
Recall that a variety $X$ is $\QQ$-Gorenstein if it is Cohen-Macaulay
and there is a integer $m$, such
that the saturation $\omega_X^{[m]}$ 
is invertible. The minimal positive $m$ satisfying this is called the
{\em canonical index} of $X$. 

\begin{definition}
Let $X$ be a normal quasi-projective $\QQ$-Gorenstein variety.
The {\em canonical covering stack}  $\cX\to X$ is defined as the
stack-theoretic quotient  
$$\cX \ =\ \  [\,P\,/\, \GG_{m}\,],$$ where  $$P\   =\ \ 
 {Spec\,}_{{Y}}\left(\bigoplus_{i \in Z} \omega_{X}^{[i]}\right).$$  
\end{definition}

The variety $P$ is the canonical $\GG_m$-space of $X$, and is
Gorenstein.
The canonical covering stack is  Gorenstein, and automatically a
Deligne--Mumford 
stack, since the stabilizer of any point of $P$ is contained in the
finite group $\bmu_m$.

More traditionally, a different but equivalent construction has been used.
For each
point $x \in X$, one can find an open neighborhood $U_{x}$ such that
$m_{x}K_{X}$ is a Cartier divisor for a minimum positive integer
$m_{x}$. The canonical covering $\pi_{x}: \tilde{U}_{x} \rightarrow U_{x}$ is
a finite morphism of degree $m_{x}$ where $\tilde{U}_{x}$ is a normal variety
and $\pi_{x}$ is \'etale in $\codim= 1 $  such that 
$K_{\tilde{U}_{x}}=\pi_{x}^{*}(m_{x}K_{X})$ is Cartier.
The canonical coverings are \'etale locally unique. Thus one can
define the canonical covering stack $\cX$ by 
the atlas given by the collection of canonical coverings. 


Consider the natural morphism $\pi:\cX \rightarrow X$.
The  stacky points of $\cX$ are the preimage of   points where $K_{X}$
 fails to be Cartier. When $X$ is a terminal
 $\QQ$-Gorennstein variety, 
the Weil divisor $K_{X}$ is Cartier in codimension $2$; in particular,
 when $\dim X =3$ this locus consists of isolated points. 

\subsection{Perverse moduli stack for $\QQ$-Gorenstein flops}
\label{Sec:thankbr}  
Let $X$ and $Y$ be terminal $\QQ$-Gorenstein threefolds and
 $f: X \rightarrow Y$ be a flopping contraction. 
Under these assumptions, the varieties $P_{X}$ and $P_{Y}$ are
 Gorenstein fourfolds. There is a    
 $\CC^{*}$-equivariant lifting $f_{P_{Y}}: P_{X} \rightarrow
 P_{Y}$. This morphism satisfies the conditions (B.1) and 
(B.2). 
Therefore we can define the moduli space of perverse point
sheaves as in \cite{br}. There is a distinguished component
 $W(P_{X}/P_{Y})$, which is birational to
$P_{Y}$. By its universal property, the action of $\CC^*$ on $P_X$ and
 $P_Y$ induces an action on 
 $W(P_{X}/P_{Y})$. 

We note that $W(P_{X}/P_{Y})$ is the flop of $P_X \to P_Y$. This follows from
the results of \cite{Chen}, as follows: let $H_Y$ be a general hyperplane
section of 
$P_Y$, with inverse image $H_X \subset P_X$. The morphism $H_X \to H_Y$ is a
Gorenstein flopping contraction. The restriction of $W(P_{X}/P_{Y})$ to 
$H_Y$ is isomorphic to $W(H_X/H_Y)$ by Proposition 4.4 of \cite{Chen}. By
Theorem 1.1 of \cite{Chen} the
latter is the flop of $H_X \to H_Y$, which gives the claim.

  Consider the stack $$ \cW \ := \ \ [\,W(P_{X}/P_{Y})\ /\
\CC^{*}\,].$$ 
By definition, we can give an interpretation of this as a moduli stack. 
For each scheme $S$, the 
 quotient stack $[W(P_{X}/P_{Y})/\CC^{*}]$ gives a category
 $\cW(S) =\ [\,W(P_{X}/P_{Y})\,/\,\CC^{*}\,](S)$.  An object in this category
consists 
of a 
$\CC^{*}$-bundle $Q \to S$ with 
 a $\CC^{*}$ equivariant morphism $\alpha: Q \rightarrow
 W(P_{X}/P_{Y})$. For any two objects $(Q_{1}, \alpha_{1})$ and  $(Q_{2},
 \alpha_{2})$  in  $\cW(S)$, an arrow is a
  $\CC^{*}$-equivariant morphism $\beta: Q_{1} \rightarrow Q_{2}$
 such that $\alpha_{1}= \alpha_{2} \circ \beta$.  

Using  the definition of
 $W(P_{X}/P_{Y})$, we can get a more concrete interpretation of
objects. 
Since $W(P_{X}/P_{Y})$ is a fine moduli space, an element in
$W(P_{X}/P_{Y})(Q)$ is equivalent to a family of perverse 
point sheaves for $P_{X} \rightarrow P_{Y}$ 
     parametrized by $Q$. Thus   
    an object in the category $\cW(S)$ is a pair
 $(Q,E)$ where $Q$ is a $\CC^{*}$-bundle over $S$ 
and $E$  is a $\CC^*$-equivariant  family of perverse point sheaves
 for $P_{X} \rightarrow 
P_{Y}$ parametrized by $Q$.  

This interpretation relies on the particular presentation of $\cW$ we
chose. A similar interpretation in terms of perverse point sheaves on
a more egeneral presentation with descent data can be obtained.

By \cite{Chen}, Proposition 4.2, the universal perverse point sheaf
over $P_X \times W(P_X/P_Y)$ is the 
structure sheaf of the fibered product $P_X \times_{P_Y}
W(P_X/P_Y)$. The natural $\CC^*$ action defines a corresponding sheaf $\cE$
on the quotient stack $\cW$, which is evidently the structure sheaf
 of
the fibered product $\cX \times_\cY \cW$. 

Can we view $\cW$ as a moduli stack of perverse point sheaves on
$\cX/\cY$? The answer is ``yes'' if we are careful about the
definition, and still a bit unsatisfactory. 
The fibers of $\cE$ over geometric points of $\cW$ are indeed
elements of $Per(\cX/\cY)$, defined just as in \cite{br}. However
the universal perverse point sheaf is {\em not} a quotient of
$\cO_{\cX \times \cW}$, for the same reason that the universal
non-perverse ``point sheaf'' $\cO_\cY$ of a non-representable stack
$\cY$ is not a quotient of $\cO_{\cY
\times \cY}$.   Only after pulling back to an \'etale cover of $\cY$
can we write it as a quotient. (One might suspect that any moduli construction
of ``Hilbert scheme'' type on a stack should requires a nontrivial game with
descent data, as we presented above, in contrast with the ``Quot
scheme'' construction of 
\cite{Olsson-Starr}.) 
 
\subsection{Proof of Proposition \ref{Prop:flop}} 
We claim that $\cW$ is the flop of $\cX \to \cY$ and the sheaf $\cE$
gives an equivalence of derived categories.

 Consider a  scheme $Y'$ and an
 \'etale morphism $Y' \to \cY$. Write $X' = \cX \times_\cY Y'$.  The
 scheme $Y'$  is a 
 terminal Gorenstein 
 threefold, and $X'\to Y'$ is a flopping contraction. Write $W' = W(X'/Y')$.
We have
\begin{lemma}
  $W' = \cW \times_\cY
 Y'$ 
\end{lemma}
{\bf Proof.} 
The formation of the category $Per$ commutes with flat base change (see e.g.,
\cite{Vandenbergh}, Proposition 3.1.4). Consider the fiber diagram
$$\xymatrix{
P_X \times_{\cX}X'\ar[rrr]\ar[rdd]\ar[d]&&&P_X\ar[rdd]\ar[d] \\
P_Y \times_{\cY}Y'\ar[rrr]\ar[rdd]      &&&P_Y\ar[rdd] \\
&X'\ar[rrr]                       \ar[d]&&&\cX        \ar[d]\\
&Y'\ar[rrr]                             &&&\cY        
}
$$ where all the horizontal arrows are flat, and the diagonal arrows are
$\GG_m$-bundles. We therefore have morphisms  
$$W(P_X/P_Y)\ \times_\cY\ Y'\ \  \lrar \ \ 
 W(\, P_X\times_\cY Y'\ /\ P_Y\times_\cY Y'\,)\ \  \longleftarrow \ 
 W'\ \times_{Y'} \  ( P_Y\times_\cY Y').$$ These are isomorphisms since all
these varieties  are the relevant flops. Also $$   W(\, P_X\times_\cY Y'\ /\
P_Y\times_\cY  Y'\,) \to W'$$ is a $\CC^*$-bundle.
Taking quotients by $\CC^*$
gives the result. \qed

 By definition the sheaf $\cE$ pulls back to the structure
 sheaf of $X'\times_{Y'} W'$, which is the universal perverse point
 sheaf for $X'/Y'$. 

Applying \cite{Chen}, Threorem 1.1, we have
 that the Fourier--Mukai type transform $D^b(W')\to D^b(X')$ is an
 equivalence and $W'\simeq {X'}^+$, the flop of $X' \to Y'$. Applying
 \cite{Chen}, Proposition  3.2 we have that the
 transform $$F_\cE\colon  D^b(\cW)\ \lrar\ D^b(\cX)$$ is an
 equivalence.

Let $D$ be a negative divisor for $f:X/Y$, let $\cD$ be its pullback to
$\cX$ and $D'$  the pullback to $X'$.  Since $W'\simeq {X'}^+$,
we have that $$W' \  =\ \  {\cP}roj\,_{Y'}\ \bigoplus_{n\geq 0}
f'_*\cO_{X'}(nD').$$ It follows that 
$$\cW  \  =\ \  {\cP}roj\,_{\cY}\ \bigoplus_{n\geq 0}
f'_*\cO_{\cX}(n\cD).$$ This completes the proof. \qed


\section{Threefold Francia flips}  
\subsection{Local Models}\label{localmodel}\label{Sec:local-models}
In this section, we study the Francia flip case.
Recall the simplest sequence of flips (see \cite{km98}, p.39).
Let  $Y_{1}$ be the quadratic singularity $\{ xy=uv \} \in \CC^{4}$. 
Denote by $X_{1}$ the variety obtained by blowing up the ideal $(x,v)$. 
Denote by $X_{1}^{+}$ the  variety obtained by    blowing up the ideal $(x,u)$.
 Let   $\Zn$ be the cyclic group of $n$ elements,  acting on $Y_{1}$
via $(x,y,u,v) \mapsto (\zeta x , y, \zeta u , v)$, which lifts to an
action on $X_1$ and on $X_1^+$. The corresponding quotients are denoted by $X_{n}$, $Y_{n}$, and $X^{+}_{n}$. 
Consider the diagram 
$$\begin{array}{ccccc} X_n & & & & X_n^+ \\
 & \searrow && \swarrow & \\
&& Y_n. &&
\end{array}$$
Note that the Picard numbers satisfy $\rho(X_{n}^{+}/Y_{n})=
\rho(X_{n}/Y_{n})=1$  
and the variety $X_{n}^{+}$ is smooth.
Denote by $C_{n}$ the exceptional curve in $X_{n}$ and by $C'_{n}$ the exceptional curve in $X^{+}_{n}$.
A standard computation shows that \[C_{n} \cdot K_{X_{n}}=
\frac{-(n-1)}{n},\] and \[C'_{n} \cdot K_{X_{n}^{+}}= n-1.\] Thus
$X_{n}^{+} \rightarrow Y_{n}$ is the flip of $X_{n} \rightarrow Y_{n}$
when $n \geq 2$.
 
\begin{definition}\label{localfrancia}
Let $f: X \rightarrow Y$ be a flipping contraction between two
quasi-projective threefolds. We call $f:  X \rightarrow Y$ a local Francia
flipping contraction if $f$ is \'etale locally isomorphic to some $X_{n}
\rightarrow Y_{n}$. 
\end{definition}
By ``\'etale locally isomorphic to  $X_{n}
\rightarrow Y_{n}$'' we mean that every point $y\in Y$ has an \'etale
neighborhood $\phi:U \to Y$  with another \'etale morphism $\psi:U \to Y_{n}$
and an isomorphism
$$X \times_Y U \simeq X_n\times_{Y_n} U.$$ This induces an isomorphism 
$$X^+ \times_Y U \simeq X^+_n\times_{Y_n} U.$$

\'Etale locally we have the following diagram 
$$
\xymatrix{ X_{1} \ar[rrd] \ar[d]&& & & X_{1}^{+} \ar[lld] \ar[d] 
\\ 
         [\,X_{1}\,/\,(\Zn)\,] \ar[rrd] \ar[d] && Y_{1}  \ar[d]  &&
[\,X_{1}^{+}\,/\,(\Zn)\,] \ar[lld] \ar[d]  \\ 
         \cX_{n}  \ar[rd]&& [\,Y_{1}\,/\,(\Zn)\,]  \ar[d] & & X_{n}^{+} \;
.\ar[lld] \\ 
 &X_{n} \ar[r] &  Y_{n}
} 
$$

Dropping the subscripts $n$ and concentrating on the left side, we have

\begin{equation} \label{Eq:quotient-diagram}
\xymatrix{&&
X_{1} \ar[lld]_{q}\ar[lldd]^{\tau}\ar[ddd]^{q'}\\  [\,X_{1}\ /\ (\Zn)\,]
\ar[d]_{\sigma}\\ \cX \ar[drr]^{\pi}\\&& X 
\ar[d]^{\bar f}\\&& Y.} 
\end{equation}

The quotient morphism   $$q\colon \ X_{1} \ \ \rightarrow\ \  [\,X_{1}\ /\
(\Zn)\,]$$ is, by definition, 
\'etale. The morphism 
$$\sigma\colon \ [\,X_{1}\ /\ (\Zn)\,]\ \  \rightarrow \ \ \cX$$ is
the natural morphism 
``forgetting the stacky structure'' along a $\Zn$  divisor
$[\,D\,/\,(\Zn)\,]$. The divisor $D\subset X_1$ is given, in the local
coordinates 
introduced 
earlier, by the function $\frac{x}{v}$ in the affine chart $$V\ \ =\
\  \Spec\  
\CC\left[\frac{x}{v},y,v\right].$$     We denote by $D'$ the image of
$D$ in 
$X$, so ${q'}^*D' = nD$. 

The exceptional curve in $X_1$ is given in this chart by $y=v=0$.

Note that the morphism $\sigma$ is  flat since $[X_{1}/ (\Zn)]
\rightarrow \cX$ is surjective, proper and quasi-finite, $\cX$ is smooth and
$[X_{1}/(\Zn)]$ is Cohen-Macaulay (in fact, it is smooth). The coarse moduli
space of $[X_{1}/(\Zn)]$ is also $X$. We also have that $\sigma_{*} \sigma^{*}
=id$ on any quasi-coherent sheaf, since $\sigma_* \cO_{  [X_{1}/ (\Zn)]} =
\cO_\cX$.

Consider now the other affine chart $$U\ \ =\ \  \Spec\
\CC\left[x,u,\frac{v}{x}\right]$$ on $X_1$. The exceptional curve is defined by
$x=u=0$. Denote by $E$ the divisor defined by 
$\frac{v}{x}$ - it meets the exceptional curve properly. Again we write $E' =
q'_* E$, so ${q'}^* E' = nE$. 

On $X_1$ we have an equality of divisors 
$$div \left(\ \frac{v}{x}\ \right)\ \ = \ \ E - D.$$ 
The function $\frac{v}{x}$ is a $\Zn$--semi-invariant  with respect to the
 character $$\frac{v}{x}\ \  \mapsto\ \
\zeta_n^{-1}\cdot\frac{v}{x}.$$ Therefore we can 
identify the $\Zn$-sheaf $L_{\chi_i}$ isomorphic to  $\cO_{X_1}$ 
twisted with action given by a character $\chi_i$ as 
$$L_{\chi_i} \mathop{\simeq}\limits^{\Zn} \cO_{X_1}(j(E-D))$$ for 
suitable $i$. 

\subsection{A perversity for Francia flips}
Let $X$ be as in Section~\ref{localmodel}. 
In this case the natural canonical covering stack $\cX$ is a smooth
Deligne-Mumford 
 stack,   whose coarse moduli space is $X$. 
 Consider the natural morphism $\pi:\cX \rightarrow X$.
The  stacky points of $\cX$ are the preimage of   points where $K_{X}$
fails to be Cartier. This loci consists of isolated points.

Let $f: \cX \rightarrow Y$ be the morphism as above. 
We write  $f = \bar f \pi $ with $\pi : \cX \rightarrow X$ and
$\bar f: X \rightarrow Y$. 

Following \cite{br},  we can define for any $m$  a perverse
$t$-structure $t(m)$  for the morphism of {\em schemes} $\bar f: X \rightarrow
Y$.  
The heart
of this $t$-structure is denoted by $\Per^{m}(X/Y)$.

  Note that $\bL \pi^{*}: D(X) \rightarrow D(\cX)$ is a
fully-faithful embedding. This follows from the projection formula, given that
$\pi_*\cO_\cX = \cO_X$.

We proceed by defining a perverse $t$-structure for the {\em stack} $\cX$.
We  define two sub categories of $D(\cX)$:
\[B=\{ \bL \pi^{*}C \in D^{b}(\cX) |\; C \in D(X)\},\] and
\[C_{2}= \{ C \in D(\cX) |\; \bR \pi_{*}C=0\}.\]

The pair $(B,C_{2})$ gives a semiorthogonal decomposition on $D(\cX)$.

On the category $C_{2}$,  the standard $t$-structure induced by the
standard $t$-structure on $D(\cX)$ gives a $t$-structure on $C_{2}$.   

Since $\bR \pi_{*}$ has the right adjoint $\pi^{!}$ and the left
adjoint $\bL \pi^{*}$, we can glue any $t$-structures on $C_{2}$ and $D(X)$.  
We define new t-structures $t(m,n)$ on $D(\cX)$: 

\begin{definition}
Denote by $t(m,n)$ the t-structure obtained by gluing:
the perverse t-structure $t(m)$ on $D(X)$,   
 and the 
standard t-structure   shifted by  $n$ on $C_{2}$.
We denote the heart of this t-structure by $\Per^{m,n}(\cX/Y)$.
\end{definition}

Inspired by \cite{br}, 
we shall only be interested in the case that $(m,n)=(-1,0)$ and denote
$\Per^{-1,0}(\cX/Y)$ simply by $\Per(\cX/Y)$.

By definition, an object $E \in  D(\cX)$ is in  $\Per(\cX/Y)$ if and only if
the following conditions are satisfied: 
\begin{enumerate}
\item $\bR \pi_{*} E$ is a perverse sheaf for the morphism $\bar f:X
  \rightarrow Y$, and
\item     \begin{enumerate} \item   $\Hom(E,C)=0$ for all $C$ in $C_{2}^{>0}$,
and 
                    \item  $\Hom(D,E)=0$ for all $D$ in $C_{2}^{<0}.$
           \end{enumerate}
 \end{enumerate}

 
\begin{definition}
A perverse sheaf $I$  is called a perverse ideal sheaf if there is an injection
$I \hookrightarrow \cO_{\cX}$ in the abelian category $\Per(\cX/Y)$. A perverse
sheaf $E$ is called a perverse structure sheaf if there is a surjection
$\cO_{\cX} \rightarrow E$ in the  
abelian category $\Per(\cX/Y)$.  
A perverse point sheaf is a perverse structure sheaf which is numerically
equivalent to the structure sheaf of a point $x \in \cX$. 
\end{definition}

\subsection{Characterization of perverse sheaves and perverse ideal sheaves}

The next lemma is an analogue of Lemma 3.2 in \cite{br} and can be proved in a
similar manner. 
\begin{lemma}\label{ps}
An object $E$ of $D(\cX)$ is a perverse  sheaf if and only if the following
four conditions are satisfied: 
\begin{enumerate}[{\rm (PS.}\rm 1{)}]
\item $H_{i}(E)=0$ unless $i=0$ or $1$,
\item $\bR^{1}f_{*}H_{0}(E)=0$ and $\bR^{0}f_{*}H_{1}(E)=0$,
\item $\Hom (\pi_{*} H_0(E),C)=0$ whenever $\bar f_* C = \RR^1 f_*C = 0$, and
\item $\Hom(D, H_{1}(E))=0$ for any sheaf $D$ in $C_{2}$ (i.e. $\pi_{*}D=0$). 
\end{enumerate} 
\end{lemma}

We continue  following \cite{br}:

\begin{lemma}\label{sheaf}
A perverse ideal sheaf is a sheaf.
\end{lemma}
{\bf Proof.} Let $F$ be a perverse ideal sheaf on $\cX$ and $E$ be the
corresponding perverse structure sheaf. Consider the exact sequence in
$\Per(\cX/Y)$ 
\[ 0 \rightarrow F \rightarrow \cO_{\cX} \rightarrow E \rightarrow 0. \]
By Lemma~\ref{ps} all the terms in this sequence have homology only in degrees
$0$ 
and $1$, and $H_1 (\cO_\cX)=0$.
Applying the homology functor to this exact sequence, we get 
$$0 \to H_1F \to 0 \to H_1(E) \cdots $$ so $H_1(F) = 0$.
\qed

\begin{lemma}\label{idealsheaf}
A sheaf $ F \in Coh(\cX)$ is a perverse ideal sheaf if and only if it satisfies
the following two conditions: 
\begin{enumerate}[{\rm (PIS.}\rm 1{)}]
\item \label{(PIS.1)} $\bR \pi_{*}F$ is a perverse ideal sheaf for $\bar f:X
  \rightarrow Y$, and
\item  \label{(PIS.3)}  $\Hom(D, F)=0$ for any sheaf $D \in C_{2}$. 

In addition, a perverse ideal sheaf satisfies  the following property:
\item \label{(PIS.2)}  the cokernel of the natural morphism $f^{*}f_{*}(F)
\rightarrow F$ 
  is in the category $C_{2}$.

\end{enumerate} 

\end{lemma}
{\bf Proof.} Let $F \in D(\cX)$ be a perverse ideal sheaf and $E$ be
the corresponding perverse point sheaf. Since $\pi_{*}$ is exact on
the abelian category $\Per(\cX/Y)$, there is an exact sequence in $Per(X/Y)$
\[ 0 \rightarrow \pi_{*}F \rightarrow \cO_{X} \rightarrow \pi_{*}(E)
\rightarrow 0. \]  
It follows that $\pi_{*}F$ is a perverse ideal sheaf for $ \bar f: X
\rightarrow Y$.
Since $\pi_{*}F$ is a perverse ideal sheaf for $\bar f: X \rightarrow
Y$, it follows that $\bar f^{*} \bar f_{*}(\pi_{*}F) \rightarrow
\pi_{*}F$ is surjective. Applying  the right-exact functor $ \pi^{*}$, we have 
$f^{*} f_{*}(F) \rightarrow \pi^{*} \pi_{*}F$ is
surjective. Therefore, to show (PIS.\ref{(PIS.2)}), it suffices to show that
the 
cokernel $D$ of   
$ \pi^{*} \pi_{*}F \rightarrow F$  is in $C_{2}$. But $ \pi_*\pi^{*} \pi_{*}F
\rightarrow \pi_*F$ is an isomorphism, therefore its cokernel $\pi_*D$
vanishes. 

To check the property (PIS.\ref{(PIS.3)}), we use the condition (PS.4) on $E$
 and 
 the exact sequence
\begin{equation}\label{Eq:homology-sequence}
 0 \rightarrow H_{1}(E) \rightarrow F \rightarrow \cO_{\cX}
\rightarrow H_{0}(E) \rightarrow 0.\end{equation}
 It follows  that $\Hom(D,F)=0$ for all sheaf $D \in C_{2}$ since we have
 $\Hom(D, H_{1}(E))=0$ and $\Hom(D, \cO_\cX)=0$.

For the converse, consider a
sheaf $F$ satisfying 
(PIS.\ref{(PIS.1)}) and 
(PIS.\ref{(PIS.3)}). Then (PS.1) is automatic as $H_i(F) = 0$ for $i\neq
0$. Since 
$\pi_*$ is exact we have $ \bR^if_* F = \bR^i\bar f_* (\pi_*F)$, and (PS.2)
follows since $\pi_*F$ is in $Per(X/Y)$. For the same reason (PS.3) holds. And
(PS.4) is exactly (PIS.\ref{(PIS.3)}). Thus $F$ is a perverse sheaf.

We now show that $\Hom(F, \cO) = \Hom(\pi_*F, \cO)$, hence
nonzero. Consider the exact 
sequence of sheaves $$0 \to A \to \pi^*\pi_*F \to F \to D \to 0,$$ 
and denote the image of $\eta: \pi^*\pi_*F \to F$ by $C$. Since $\eta$ is an
isomorphism away from the singular locus of $X$, the sheaf $A$ is torsion,
therefore $\Hom(C, \cO) \subset \Hom (\pi^*\pi_*F, \cO)$. Also the sheaf $D$ is
torsion, and we have an exact sequence
$$0 \to \Hom(F, \cO) \to \Hom(C, \cO) \to \Ext^1(D, \cO),$$
but since $D$ is supported in dimension 0, $$\Ext^1(D, \cO) = H^2(D \otimes
\omega_\cX)^\vee = 0.$$    

Fix a nonzero element in  $\Hom(F, \cO)$ and consider the triangle 
$F \to \cO \to E \to F[1]$. It suffices to show that $E$ is perverse. The long
exact 
sequence of homology (sequence \ref{Eq:homology-sequence} above)
 gives (PS.1). Clearly
$\bR\pi_* E$ is the perverse quotient of the corresponding element of
$\Hom(\pi_*F, \cO)$, so (PS.2) and (PS.3) follow. And (PS.4) follows again
because $\Hom(D, H_1(E)) \subset \Hom(D, F) = 0$ by  (PIS.\ref{(PIS.3)}).
 \qed

Since $\cX \rightarrow Y$ is an isomorphism outside the singular
points of $Y$, we can compactify $\cX$ and $Y$ to get  $\bar f: \bar \cX
\rightarrow 
\bar Y$. 
It is clear from the definition  that perverse point sheaves are local
objects over $\bar Y$, i.e. $W(\bar \cX /\bar Y)_{Y} \cong W(\cX/Y)$. Abusing
the 
notation, we shall still denote this new morphism by 
 $f: \cX \rightarrow Y$. 

\subsection{Simplicity of perverse point and perverse point-ideal  sheaves}

\begin{lemma}\label{pointsheaf}
Let $E_{1}$ and $E_{2}$ be two perverse point sheaves. Then,

\[ \Hom(E_{1}, E_{2})= \left\{ \begin{array}{r@{\quad:\quad}l} 
      0 & \text{if} \; E_{1} \not\cong E_{2}, \\  
\CC& \text{if} \; E_{1} \cong E_{2}. \end{array}\right. \]
\end{lemma}

{\bf Proof.}  See Lemma 3.6 in \cite{br}, using the following:

\begin{lemma}\label{numericalto01}
Let $E \in D^{b}(\cX/Y)$ be a perverse sheaf with $Supp(E) \subset
f^{-1}(y)$ for a point $y \in Y$. If $E$ is numerically equivalent to
$0$, then $E \cong 0$.

\end{lemma}
\textbf{Proof.} This follows immediately from the following
Lemma~\ref{numericalto00} (which is implicit in \cite{br}) and property  
(PS.4) in Lemma~\ref{ps}. $\Box$ 

\begin{lemma}\label{numericalto00}
Let $E \in D^{b}(X/Y)$ be a perverse sheaf with $Supp(E) \subset
\bar f^{-1}(y)$ for a point $y \in Y$. If $E$ is numerically equivalent to
$0$, then $E \cong 0$.
\end{lemma}

\textbf{Proof.}   Let $E$ be a perverse sheaf which is numerically
equivalent to $0$. It suffices to show that all homology groups vanish, that
 is $H_{0}(E)=H_{1}(E)=0$. 

First note that $\chi(H_{0}(E))= \chi(H_{1}(E))$. Since $E$
is a perverse sheaf, we have $\bR^{1} \bar f_{*} (H_{0}(E))= 
\bar f_{*}(H_{1}(E))=0$. Therefore, $\chi(H_{0}(E))$ is a
nonnegative integer and $\chi(H_{1}(E))$ is a nonpositive
integer. This implies   $\chi(H_{0}(E))= \chi(H_{1}(E))= 0$. By the support
assumption, $f_{*}(H_{0}(E))$ is supported at a point, and this
implies that $H^0 \bar f_{*}(H_{0}(E))=0$, which means $H_{0}(E)=0$. The same
holds for $\bR^{1} \bar f_{*} (H_{1}(E))$. 

We thus  have 
$\bR \bar f_{*}(H_{i}(E))=0$ and therefore the two sheaves $H_{i}(E)$ have
support in 
pure dimension 1. 
For $H_0(E)$ we have  $\Hom(H_{0}(E),H_0(E))=0$ by (PS.3) in Lemma 3.2 in
\cite{br}. This 
implies 
$H_{0}(E) \cong 0$. 
Therefore $H_1(E)$ is a numerically trivial sheaf.

Now consider  $L$ a sufficiently ample bundle, Then$H^1 (L \otimes
H_{1}(E)) = 0 $ and $L \otimes
H_{1}(E)$ is generated by global sections, therefore $\chi( L \otimes
H_{1}(E))>0$, so $H_1(E)$ is not numerically trivial, giving a contradiction.
\qed

\begin{lemma}\label{idealsimple}
\begin{enumerate}
\item Let $F_{1}$ be a perverse point-ideal sheaf. Then
$\dim \Hom(F_{1}, \cO_\cX)=1.$
\item
Let $F_{1}$ and $F_{2}$ be perverse point-ideal sheaves. Then
\[ \Hom(F_{1}, F_{2})= \left\{ \begin{array}{r@{\quad:\quad}l} 
      0 & \text{if} \; F_{1} \not\cong F_{2}, \\  
\CC& \text{if} \; F_{1} \cong F_{2}. \end{array}\right. \]
\end{enumerate}
\end{lemma}

{\bf Proof.} Let $F_{1}$ be a perverse  point-ideal sheaf. 
Consider the following exact sequence (in the usual abelian category)
\[ 0 \rightarrow  A \rightarrow f^{*}f_{*}F_{1} \rightarrow F_{1} \rightarrow C
\rightarrow 0 \] 
where $C$ is an object in $C_{2}$ since $\pi_{*} F_{1}$ is a perverse
point-ideal sheaf for $\bar f: X \rightarrow Y$. It is also clear that $C$ is a
torsion sheaf. 
Taking $\Hom(-, \cO_{\cX})$, we get 
\[ 0 \rightarrow \Hom(F_{1}, \cO_{\cX})  \rightarrow \Hom(f^{*}f_{*}F_{1},
\cO_{\cX}). \] 
Since $$\dim
(\Hom(f^{*}f_{*}F_{1}, \cO_{\cX}))=\dim (\Hom(f_{*}F_{1}, f_{*}\cO_{\cX})) =
1,$$ this shows that $\dim (\Hom(F_{1}, \cO_{\cX}))=1$.  

Note that by (PIS.\ref{(PIS.3)}) we have  $\Hom(C, F_{2})=0$ for any sheaf  $C \in
C_{2}$. 
Taking $\Hom(-, F_{2})$, we get 
\[ 0 \rightarrow \Hom(F_{1}, F_{2}) \rightarrow \Hom(f^{*}f_{*}(F_{1}),
F_{2}). \] 
Since $\Hom(f^{*}f_{*}(F_{1}), F_{2})= \Hom(f_{*}(F_{1}),
f_{*}(F_{2}))$,  which has dimension $\leq 1$,  it follows that
$\dim (\Hom(F_{1},F_{2})) \leq 1$,  and if the dimension is $1$ then
the map factors $\Hom(F_{1}, \cO)$. Since we have $F_{1}
\rightarrow \cO$ is an injection in $\Per(\cX/Y)$, it follows that
$\theta: F_{1} \rightarrow F_{2}$ is also an injection in
$\Per(\cX/Y)$. The cokernel of $\theta$ is in $\Per(\cX/Y)$ and
numerically equivalent to $0$. Therefore it is isomorphic to $0$ by
Lemma~\ref{numericalto01}. $\Box$

\subsection{Perversity and the dualizing sheaf}

\begin{lemma}\label{cycliccover}
If $B$ is a sheaf on $\cX$ such that $B[1]$ is a perverse sheaf, then 
$f_{*}(B \otimes \omega)=0$
\end{lemma}

{\bf Remark.} The idea is that $\omega_\cX$
is negative along the exceptional 
curve, so tensoring by $\omega_\cX$ should not give more sections. This is not
quite correct as it is - it fails for torsion sheaves supported at the
non-representable point, so we need to be a bit careful and use the structure
of the sheaf $B$.

{\bf Proof.}  We have the propoerties:
\begin{itemize} 
\item $f_{*}B=0$, 
\item $B$ has pure support on the exceptional loci
\end{itemize} This in particular implies that the problem is 
\'etale local over $Y$. Thus we may 
assume that we have a diagram as in Diagram \ref{Eq:quotient-diagram}, Section
\ref{Sec:local-models}. 

For any $\Zn$-sheaf $F$ on $X_1$ we have an isotypical  decomposition 
$$q'_* F = \bigoplus_{\chi_i\colon\ \Zn\ \to\ \GG_m} (q'_*F)^{\chi_i}.$$
If $F$ is supported along the exceptional locus then there is an
isomorphism
$F\subset F(iE-iD)$ given locally by multiplication by a power of
$\frac{v}{x}$. In particular, for $i>0$ we get a $\chi_i$-twisted embedding
$F\subset F(nE-D)$.  Thus,  for any {\em nontrivial} character $\chi_i$ we can
write $$(q'_*F)^{\chi_i} \subset (q'_*F(nE\ -\ D))^\Zn.$$  

We now analyze the sheaf $q'_* \tau^* (B\otimes \omega_{\cX_1})$. First
notice that $$\tau^* (B\otimes \omega_{\cX_1}) \subset \tau^*B \otimes
\omega_{X_1}(-(n-1)D) \subset \tau^*B \otimes
\omega_{X_1},$$ which is isomorphic to $\tau^*B$ as a sheaf, but with a
different action given by twisting with a character $\chi$. This twisting only
has the effect of permuting the $\chi_i$-isotypical components of the direct
image by $q'$. 

We thus get
$$q'_* \tau^* (B\otimes \omega_{\cX_1}) = q'_* (\tau^*B(-(n-1)D)).$$
Note also that $nE - nD$ is the pullback of the principal divisor $E'-D'$ on
$X_n$ 

Applying the discussion above we get 
\begin{eqnarray*}
q'_* (\tau^*B(-(n-1)D))&\hookrightarrow& q'_* (\tau^*B)^{\Zn}\ \ \oplus\ \ 
\bigoplus_{\chi\neq 
1}\  \Big(\ q'_* \big(\,\tau^*B\,(\,-(n-1)D\ +\ nE - D\,)\,\big)\ \Big)^{\Zn} \\
&=& q'_* (\tau^*B)^{\Zn}\ \ \oplus\ \  \bigoplus_{\chi\neq
1}\  \Big(\ q'_* \big(\,\tau^*B( nE - nD)\,\big)\ \Big)^{\Zn} \\
&\simeq&  \bigoplus_{\chi:\Zn\to \GG_m}\ \ \Big(q'_* (\tau^*B)\Big)^{\Zn} \\
&\simeq& (\pi_* B)^{\oplus n}
\end{eqnarray*}

In particular we get an embedding  
$$\pi_*(B\otimes \omega) \hookrightarrow (\pi_* B)^{\oplus n}.$$

Applying the left exact functor $\bar f_*$ we get 
$$f_*(B\otimes \omega) \hookrightarrow (f_* B)^{\oplus n}  = 0.$$
\qed
  
\begin{lemma}\label{leftexact}

The functor $E \mapsto E \otimes \omega$ is perverse-left exact, i.e. sends 
$\Per(\cX/Y)$ to $D^{\geq 0}(\cX)$.
\end{lemma}
We use the next lemma to prove Lemma~\ref{leftexact}.

\begin{lemma}\label{sublemma}
Suppose $B$ is a sheaf on $\cX$ such that $B[1]$ is perverse. Then $\wB[1]$ is
also perverse. 

\end{lemma}

{\bf Proof.} Let $\pi: \cX \rightarrow X$ and $f: \cX \rightarrow Y$ be as
above. First note that $f_{*}B=0$ and $\Hom(C,B)=0$ for $C$ a sheaf in $C_{2}$,
so $B$ contains no sky-scraper sheaves. Hence the same is true for $\wB$.   

We use Lemma~\ref{ps}. The only non-trivial conditions  are the second part of
(PS.2) and (PS.4). 
The second part of (PS.2) follows from Lemma~\ref{cycliccover}.
For (PS.4), take an object $D$ such that $\pi_{*} D=0$. Since $X$ has only
isolated singularities, any such $D$ must be supported in dimension $0$.  
Note that $H_{1}( \wB[1])= \wB$.
Thus we have $\Hom(D, \wB)=0$ since $\wB$ contains no sky-scraper sheaves.

{\bf Proof of Lemma~\ref{leftexact}.} This is an easy consequence of the above
 lemma. 
 Assume the contrary. We can take a perverse sheaf $E$ and another
 object $D$ such that $\Hom(D[1], \wE)$ is nonzero. Taking homology in
 the standard t-structure, there must be a non-zero map $H_{0}(D)
 \rightarrow H_{1}(E) \otimes \omega =C$.  Since $H_{1}(E)[1]$ is a
 perverse sheaf, it follows that $C[1]$ is also a perverse sheaf by
 Lemma~\ref{sublemma}. The sheaf $H_{0}(D)$ is also a perverse sheaf
 since $D$ is a perverse sheaf. This then gives a homomorphism in
 $\Hom(H_{0}(D), C)$, a contradiction. $\Box$ \\


\begin{lemma}\label{Lem:homs}
Let $E_{i}$, $i=1,2$,  be perverse point sheaves. Then
$\Hom^{i}(E_{1},E_{2})=0$ unless $0 \leq i \leq 3$. 
\end{lemma}
{\bf Proof.}

Since
$E_i$ are objects in the abelian category $Per(\cX/Y)$, we have that
$\Hom^{i}(E_{1},E_{2})=0$ 
for $i < 0$.  

 If $\bR f_{*}(E_{1}) \neq \bR f_{*}(E_{2})$, then
$\Hom^{i}(E_{1},E_{2})=0$ for all $i$ since their supports are disjoint, so we
only need to consider the case  $\bR f_{*}(E_{1}) = \bR f_{*}(E_{2})$, and thus
the problem is local over $Y$. 

Replacing $Y$ by an \'etale base change we may assume that $\cX$ is a global
quotient stack, and thus $E_i$ are quasi isomorphic to bounded complexes of
vector bundles. Serre duality gives 
$$\Hom^{i}(E_{1},E_{2}) = \Hom^{3-i}(E_2,E_1\otimes \omega)^\vee.$$ 
This gives the result. \qed


\section{Moduli of perverse point sheaves for the Francia flips}
\subsection{The moduli space}
A family of perverse point sheaves for $\cX/Y$ parametrized by a scheme $S$ is
an 
object $\cE$ of $D^b(S \times \cX)$ such that for every point $i:s\in S$ the
fiber $Li^*\cE$ is a perverse point sheaf. Two such families $\cE, \cE'$ are
equivalent if there is a line bundle $M$ on $S$ and an isomorphism $\cE \simeq
\cE'\otimes M$. 

Define a functor
$$\cM(\cY/Y)\ \colon \ \ \cS ch\  \lrar \ \cS ets$$
which to a scheme $S$ assignes the set of equivalence classes of
families of perverse point sheaves parametrized by $S$. 

Since every perverse point sheaf determines and is determined by a
perverse point-ideal sheaf, the we can view the functor $\cM(\cX/Y)$
as the moduli functor of equivalence classes of perverse point-ideal
sheaves. By Lemma~\ref{idealsimple} the endomorphism group of such a
sheaf is $\CC$, and the automorphism group is therefore $\CC^*$.  It 
follows by standard representability theory (see, e.g., \cite{Artin}) that the
\'etale sheaf associated to $\cM(\cX/Y)$ is represented by an algebraic space
$M(\cX/Y)$, locally of finite type over $\CC$. An argument of 
Mukai 
(see, e.g. \cite{br}, proof of Theorem 5.5) shows that  $M(\cX/Y)$ is
 a fine moduli space for $\cM(\cX/Y)$, i.e. there is a universal perverse
point-iedal sheaf over $M(\cX/Y) \times \cX$.
 

\begin{lemma}
The algebraic space $M(\cX/Y)$ is separated.
\end{lemma}
\textbf{Proof.} We use the valuative criterion for separation. Let $R$
be a discrete valuation ring, with  fraction field  $K$. 
Fix a perverse ideal sheaf $F$ over $ Spec(K)$.  
Let $F_{1}$ and $F_{2}$ be two extensions of $F$  to $Spec(R)$ with an
isomorphism $s:F_{1}|_{Spec(K)} \cong F_{2}|_{Spec(K)}$. We can
write it as $s=u^{n} \cdot h$, where $u$ is a uniformizer in $R$ and
$h$ extends comes from  a homomorphism $F_1 \to F_2$. Taking $n$
minimal, we may assume that the restriction of $\bar h: F_{1\,s} \to
F_{2\,s}$ of  $h \in 
\Hom(F_{1},F_{2})$, to the special fiber is nonzero. By Lemma
\ref{idealsimple} we have that $F_{1} 
\cong F_{2}$ over $Spec(R)$. \qed  
 
\subsection{Projectivity}

There is a
 distinguished component of $M(\cX/Y)$, which is birational to $Y$. We
 shall denote this component by $W$.
To complete Theorem \ref{Th:francia-W} we need to show that the
 distinguished component $W$  is projective over $Y$.

  Note that it suffices to consider the case that $Y \cong
Y_{n}$ and $X \cong X_{n}$.  
 By Lemma \ref{Lem:flip-bir} below, there is a birational morphism
$X_{n}^{+} 
\rightarrow W$; and  the composition
$X_{n}^{+} \rightarrow W \to  W(X_{n}/Y_{n})$ is a
finite morphism, where $W(X_{n}/Y_{n})$ is the moduli space for the usual
perverse point sheaves. This  implies  that     $W \rightarrow
W(X_{n}/Y_{n})$  
is finite. Since $ W(X_{n}/Y_{n}) \rightarrow Y_{n}$ is projective, it 
follows immediately that $W \rightarrow Y_{n}$ is also projective. \qed

\begin{lemma} \label{Lem:flip-bir} There exists a  birational morphism
$X^+ \to W$. The composition $X^+ \to W \to W(X/Y)$ is finite.
\end{lemma}

{\bf Proof.}
To give  a morphism $X_{n}^{+} \rightarrow W$,  
 we  exhibit a family of perverse point sheaves
over $X_{n}^{+}$.
 Unlike the terminal
Gorenstein case, the correct family of perverse point sheaves  is not the fiber
 product $X_{n}^{+} 
\times_{ Y_{n}} \cX$. This fiber product has an extra embedded component. We
shall show that the correct candidate is the structure sheaf of the
reduction of the fiber 
product. We denote  
it by 
 $\widetilde{(X^{+}_{n} \times_{Y_{n}} \cX_{n})}.$

To show this is a
 family of perverse point sheaves, we use Lemma~\ref{idealsheaf} and check the
 conditions (PIS.\ref{(PIS.1)})-(PIS.\ref{(PIS.3)}) for the corresponding
 perverse ideal sheaves. 

Consider  the morphism $id\times \pi: X^{+}_{n}
\times_{Y_{n}} \cX_{n}) \to X_n^+ \times_{Y_{n}} X_n$. Condition
(PIS.\ref{(PIS.1)}) amounts to checking that $(id\times \pi)_*
I_{\widetilde{(X^{+}_{n} \times_{Y_{n}} \cX_{n})}}$ is a perverse ideal
sheaf. Since $\pi$ is exact, this is nothing by the ideal sheaf
$I_{\widetilde{(X^{+}_{n} \times_{Y_{n}} X_{n})}}$ of the fibered product. In
\cite{A-C}, Proposition 1.3.2, it is shown that this is indeed a perverse ideal
sheaf.

To check the condition (PIS.\ref{(PIS.3)}), we consider  
 the flat base change $p:X_{1} \rightarrow \cX$. We have a natural 
morphism
$p_{n}:X_{1} \times_{Y_{1}} X_{1}^{+} \rightarrow X_{1} \times X_{n}^{+}$. 
Denote by $Z$ the image of this morphism in $X_{1} \times X_{n}^{+}$. 
The  morphism $p_n$
 factors through 
\[
X_{1} \times_{Y_{1}} X^{+}_{1} \rightarrow
 X_{1} \times_{Y_{n}} X_{n}^{+} \hookrightarrow  X_{1} \times X_{n}^{+}.\]
The first morphism is finite and the second morphism is a closed
embedding. Therefore the composite  is a  finite morphism. Note also that
this morphism is birational, since the projections to $X_{1}$ are
birational. 
\begin{lemma}
The morphism  $p_{n}:X_{1} \times_{Y_{1}} X_{1}^{+} \to Z $ is an isomorphism.
\end{lemma}

{\bf Proof.}
It suffices to check that the image $Z$ in $X_{1} \times X_{n}^{+}$ is
normal by Zariski's main theorem.
Indeed, we shall prove  that  $Z$ is  smooth.
To check this, we do  explicit computations on affine charts.

The variety $Y_{1}$ is given by  $\{xy -uv=0\} \subset \CC^{4} $. The
 variety
 $X^{+}_{1}$ is obtained by blowing up
the ideal $(x, u)$, and the variety $X_{1}$ is obtain by blowing up
$(x,v)$. The variety $X_{1}^{+}$ can be cover by two affine pieces $U^{+}_{u
 \not= 0}$ 
 and $U^{+}_{x \not= 0}$. Similarly, the variety $X_{1}$ can be cover by two
 affine pieces $U^{-}_{v \not= 0}$ 
 and $U^{-}_{x \not= 0}$. 
The variety $X_{1} \times X_{1}^{+}$ can be covered by four
affine charts. 
 
We first consider the affine chart $U^{+}_{u \not= 0}$ in $X_{1}^{+}$
and $U^{-}_{v \not= 0}$ in $X_{1}$.
The structure sheaf $\cO_{Y_{n}}$ is generated by $\{z_{i}=x^{i}u^{n-i}, y, v
\}$.   
The structure sheaf $\cO_{1}$ is generated by
$\{y_{-},v_{-},s_{-}=x_{-}/v_{-}\}$. 
The structure sheaf $\cO_{1}^{+}$ is generated by $\{t_{+}=x_{+}/u_{+}, y_{+},
u_{+}\}$. 
The structure sheaf
$\cO_{X_{n}^{+}}$ is generated by $\{t_{+}, y_{+}, u^{n}_{+}\}$. 

The ideal of the  fiber product $X_{n}^{+} \times_{Y_{n}} X_{1}$ in
$X_{n}^{+} \times X_{1}$
 is 
generated by
$\{t_{+}y_{+}-v_{-},y_{-}-y_{+},t_{+}^{i}z_{+}-s_{-}^{n}v_{-}^{i}y_{-}^{n-i}\}
$. An easy calculation shows that this ideal is generated by 
$\{t_{+}y_{+}-v_{-},y_{-}-y_{+}, z_{+}-s_{-}^{n}y_{-}^{n}\}$. The variety
defined by this ideal is smooth, and it is isomorphic to $Z$ in this affine
chart.

On the second chart $U^{+}_{u \not= 0}$ in $X_{1}^{+}$ and $U^{-}_{x \not= 0}$
in $X_{1}$. 
The ring $\cO_{X_{n}^{+}}$ is generated by $\{t^{+}, y^{+}, u^{n}_{+}\}$
as above. The ring $\cO_{1}^{-}$ is generated by 
$\{u_{-},x_{-},s_{-}=v_{-}/x_{-} \}$. 
The ideal of the fiber product is generated by $\{t_{+}y_{+}-s_{-}x_{-},
s_{-}u_{-}-y_{+}, t_{+}^{i}z_{+}-x_{-}^{i}u_{-}^{n-i}\}$.
The fiber product is not smooth, but $Z$ is a subscheme 
of $X_{1}\times _{Y_{n}}
X_{n}^{+}$ 
and it is still smooth, since we get the ideal of the image of $X_{1}
\times _{Y_{1}} X_{1}^{+}$ in $X_{1} \times X_{n}^{+}$ after adding
the element $t_{+}u_{-}-x_{-}$ to  the ideal of the fiber product $X_{1}
\times_{Y_{n}} X_{n}^{+}$. 

For the affine chart $U^{+}_{x \not = 0}$ in $X_{1}^{+}$ and $U^{-}_{v
 \not 0}$ in $X_{1}$, it is easy to check the  variety 
 is smooth. 
For the affine chart $U^{+}_{x \not = 0}$ in $X_{1}^{+}$ and $U^{-}_{x
 \not = 0}$ in $X_{1}$, one can write down a similar
generators and equations. The ideal of the fiber product is generated by
 $\{v_{+}-S_{-}x_{-}, 
t_{+}v_{+}-s_{-}u_{-}, t_{+}^{n-i}z_{+}-x_{-}^{i}u_{-}^{n-i}\}$. It is
also not smooth.
This time we need to add the element $z_{+}-x_{-}^{n}$ to get the ideal of
the image $Z$.

\begin{lemma}
The ideal sheaf $I_{\widetilde{(X^{+}_{n} \times _{Y_{n}} \cX_{n})}}$ is flat
over $\cO_{X^{+}_{n}}$.  
\end{lemma}
\textbf{Proof} We can  check this after  pullback by the flat morphism
$X_{1}
\rightarrow \cX_{n}$. 
This amounts to checking that the ideal
$I_{Z}$ is flat over $\cO_{X_{n}^{+}}$.

First note that $I_{X_{1}^{+} \times _{Y_{1}} X_{1}}$ is flat over
$\cO_{X_{1}^{+}}$.

To check $I_{Z}$ is flat over
$X_{n}^{+}$, it suffices to show that \[Tor^{\cO_{X_{n}}}_{2}(\cO_{
    Z}, M)=0\] for all $M$. Since \[Tor^{\cO_{X_{1}}}_{2}(\cO_{Z }, M)=0,\] the
desired result would 
    follow if we can show $\cO_{X_{1}^{+}}$ is flat over $\cO_{X_{n}^{+}}$.
This holds since $X_{n}^{+}$  
is smooth and $X_{1}^{+}$ is Cohen Macaulay. $\Box$

Now we check the condition (PIS.\ref{(PIS.3)}) for the ideal sheaf
$I_{\widetilde{(X^{+}_{n} \times _{Y_{n}} \cX_{n})}}$.  Consider the natural
group of  
$\Zn$  on $X_{1}$. There is a group 
action of $\Zn$ on  $X_{1} \times X_{n}^{+}$: The group
$\Zn$ acts  on $X_{1}$ via the natural action, and acts trivially on
$X_{n}^{+}$. 
Consider the morphism $\tau:X_{1} \rightarrow \cX$ (see diagram
\ref{Eq:quotient-diagram}). Let $F$ be a sheaf on 
$\cX$. 
To check  \[ \Hom(C,F)=0\] for all $C \in D(\cX)$ satisfying $\tau_*C=0$, it
suffices to 
check that $\tau^{*}F$ has no sections supported on the preimage of
stacky points of $\cX$.

Take the family of sheaves $I_{Z}$ over the variety $X_{n}^{+}$. We
have the exact sequence
\[ 0 \rightarrow I_{Z} \rightarrow \cO_{X_{1} \times X_{n}^{+}}
\rightarrow \cO_{Z} \rightarrow 0. \] 
Let $i:p \hookrightarrow X_{n}^{+}$ be a point of $X_{n}^{+}$.
Tensoring the exact sequence  with the sheaf $\cO_{p}=\cO_{X_{n}^{+}}/m_{p}$, we get the following
\[ 0 \rightarrow \Tor^{\cO_{X_{n}^{+}}}_{1}(\cO_{p}, \cO_{Z})
\rightarrow i^{*}I_{Z} \rightarrow \cO_{X_{1}} \rightarrow
\cO_{i^{*}Z} \rightarrow 0.\]
It is clear that the   image of $i^{*}I_{Z}$
in $\cO_{X_{1}}$ is torsion free on the fiber.
To show that $i^{*}I_{Z}$ has torsion with support in pure dimension $1$, 
it suffices to show that the sheaf $\Tor^{\cO_{X_{n}^{+}}}_{1}(\cO_{p},
\cO_{Z})$ 
has support in pure dimension $1$.
To achieve this,  
we use the change of ring formula 
\[ \Tor^{\cO_{X_{1}^{+}}}_{1}(\cO_{X_{1}^{+}}/m_{p}\cO_{X_{1}^{+}}, \cO_{Z})
= \Tor ^{\cO_{X_{n}^{+}}}_{1}(\cO_{X_{n}^{+}}/m_{p} \cO_{X_{n}^{+}},
\cO_{Z}). \] 
Since $Z$ is the universal perverse point sheaf for the morphism $X_{1}
\rightarrow Y_{1}$, it 
can be directly checked that the sheaf $ \Tor^{\cO_{X_{1}^{+}}}_{1}(
\cO_{X_{1}^{+}}/m_{p}\cO_{X_{1}^{+}}, \cO_{Z})$ has  support of
pure dimension $1$.   It  
follows that the sheaf $\Tor ^{\cO_{X_{n}^{+}}}_{1}(\cO_{X_{n}^{+}}/m_{p},
\cO_{Z})$ has  support of pure dimension $1$. (Alternatively, the automorphism
group of $Z \to X_1^+$ acts transitively on the fiber, and the sheaf cannot
have sections with support in a point which is not fixed.)

Finally we verify that the composite morphism $X^+ \to W(X/Y)$ is
 finite. This composition is given by the sheaf $\cO_{X^+\times_YX}$,
 which is a perverse sheaf for $X/Y$ by \cite{A-C}, Proposition
 1.3.2. The same proposition shows that $X^+ \to W(X/Y)$ is the
 normalization, which is therefore finite.
 This completes the proof of Lemma \ref{Lem:flip-bir}. \qed

\subsection{Proof of Theorem \ref{Th:francia-FM}}
Since, by Lemma \ref{Lem:flip-bir}, we  have a finite
birational morphism $X^+ \to W$, once 
we show $W$ is smooth it follows that the morphism is an isomorphism.
We prove that $W$ is smooth and the Fourier--Mukai transform is fully
faithful at the same time, following the method of \cite{bkr}.

We denote the functor induced by the universal perverse point sheaf
$$\Phi: D(W) \to D(\cX).$$ 
As in \cite{bkr}, $\Phi$ has a left adjoint $\Psi$ and the
composition $\Psi\circ \Phi$ is defined by a sheaf $Q$ on $W\times
W$, which is supported in $W\times_YW$. We consider the restriction of
$Q$ to the complement of the diagonal. The support of $Q$ in this open
set has dimension 2, since the fibers of $W\to Y$ are images of the
fibers of $X^+ \to Y$ which have dimension $\leq 1$. Also, by Lemma
\ref{Lem:homs} the homological dimension of $Q$ on the open set is
$\leq 3$. By the intersection theorem (\cite{bkr}, Corollary 5.2) it
follows that $Q$ vanishes outside the diagonal. The argument of
\cite{bkr} Section 6, step 5-6 shows that  $W$ is nonsingular and $Q$ is a line
bundle on the 
diagonal, which implies that $\Phi$ is fully faithful. The proof is
complete. \qed

\end{document}